\documentclass[11pt,a4paper]{article}
\usepackage{amsmath}
\usepackage{amssymb}
\usepackage{latexsym}

\newtheorem{theorem}{Theorem}[section]

\newtheorem{lemma}[theorem]{Lemma}
\newtheorem{proposition}[theorem]{Proposition}
\newtheorem{definition}[theorem]{Definition}
\newtheorem{corollary}[theorem]{Corollary}

\newtheorem{exmp}[theorem]{Example}
\newtheorem{exmps}[theorem]{Examples}
\newtheorem{rem}[theorem]{Remark}
\newenvironment{example}{\begin{exmp}\rm}{\end{exmp}}
\newenvironment{examples}{\begin{exmps}\rm}{\end{exmps}}
\newenvironment{remark}{\begin{rem}\rm}{\end{rem}\rm}
\newcommand{\prf}{{\em Proof}. }
\newcommand{\qed}{\hspace*{\fill}$\Box$}

\newcommand{\beeq}[1]{\begin{eqnarray}\label{#1}}
\newcommand{\eneq}{\end{eqnarray}}

\newcommand{\ka}{{\cal A}}

\newcommand{\kc}{{\cal C}}

\newcommand{\kh}{{\cal H}}
\newcommand{\kk}{{\cal K}}

\newcommand{\IC}{{\mathbb C}}

\newcommand{\IP}{{\mathbb P}}

\newcommand{\IR}{{\mathbb R}}
\newcommand{\IZ}{{\mathbb Z}}

\newcommand{\lra}{\longrightarrow}

\newcommand{\End}{{\rm End}}

\newcommand{\id}{{\rm id}}

\newcommand{\Ker}{{\rm ker}}

\newcommand{\verylongarrow}[1]{\hbox to #1{\rightarrowfill}}

\begin{document}

{\parindent0mm{\Large\bf Infinitesimal Variation of
Harmonic Forms\\
and Lefschetz Decomposition}}

\bigskip

\bigskip

{\parindent0mm{\large\bf Daniel Huybrechts\\} {\small Mathematisches
Institut\\ Universit\"at zu K\"oln\\ Weyertal 86-90\\ 50931 K\"oln,
Germany\\}

\bigskip

\bigskip

\bigskip

\bigskip

Let $X$ be a compact K\"ahler manifold and let $\ka(X)_{cl}$ denote
the space of closed forms on $X$. Each choice of a K\"ahler form
$\omega$ on $X$ defines a Hodge decomposition of $\ka(X)_{cl}$
into the subspace of $d$-exact forms and the subspace of harmonic
forms $\kh_\omega$. How does the subspace $\kh_\omega \subset\ka(X)_{cl}$
depend on $\omega$? In this paper we study the infinitesimal
variation of $\kh_\omega$.

The most interesting situation to which the results apply and which was the
original motivation for this work is the following: Let $X$ be a Calabi-Yau
manifold and $\kk^0$ be the set of all Ricci-flat K\"ahler forms. Due to Calabi
and Yau, the natural projection to cohomology defines a bijection between
$\kk^0$ and the K\"ahler cone $\kk_X\subset H^{1,1}(X,\IR)$ of all K\"ahler
classes. Thus, varying the Ricci-flat K\"ahler structure on $X$ is equivalent
to moving in the open subset $\kk_X$ of the affine space
$H^{1,1}(X,\IR)$.
Assuming Mirror Symmetry, changing the K\"ahler class in $\kk_X$ corresponds
to deforming the complex structure of the mirror partner of $X$. The space
of all complex structures on the mirror partner is, a priori, a highly
non-linear object. So, one has tried to introduce linear coordinates
on this space corresponding to the linear coordinates on $H^{1,1}(X,\IR)$
in the A-model \cite{Ran}.
However, it turns out that passing from Ricci-flat K\"ahler forms to
K\"ahler classes one loses useful and interesting information.
In \cite{Huy} we have shown that there is a relation, though still
mysterious, between the shape of $\kk^0$, the question whether the
top exterior power of harmonic $(1,1)$-forms is harmonic, and the geometry
of $X$, in particular the existence of extremal curves.
In this paper we will show that the infinitesimal variation of the
space  $\kh_\omega$ of harmonic forms is also related to the (non)-harmonicity
of the product of harmonic forms. Thm.\ \ref{Mainthm} describes the infinitesimal
variation of $\kh_\omega$ in general terms. One immediate consequence of it is
the following result, which is a special
case of Cor.\ \ref{harmonicav}.

\bigskip

{\it Let $X$ be a compact K\"ahler manifold of dimension $N$
and let $\omega$ be a K\"ahler form
on $X$. If $\alpha$ is harmonic with respect to $\omega$,
then $\alpha$ is harmonic with respect to an infinitesimal change of
$\omega$ by a closed $(1,1)$-form $v$ if and only if the product
$\alpha v\omega^{N-2}$ is harmonic
with respect to $\omega$.}

\bigskip

As $(\kh^{1,1}_\omega)_\IR$ is the tangent space of the set $\kk^0$ of all
Ricci-flat K\"ahler forms at $\omega$, the formula in Thm.\ \ref{Mainthm}
can also be used to calculate the second fundamental form of $\kk^0$
inside $\ka^{1,1}(X)_{cl}$. Fixing a scalar product on $H^{1,1}(X,\IR)$
the mean curvature then associates to each Ricci-flat K\"ahler form a
unique exact $(1,1)$-form.

The method to prove these results is almost entirely algebraic and of some
interest in itself. Due to the fortunate interplay between the Hodge
$*$-operator, needed to define the Laplacian, and the Lefschetz operator,
one can make use of general results on the variation of $sl_2$-representations
presented in Sect.\ \ref{slgeneral}. Let
\(\displaystyle {\oplus_{n=0}^{2N} A^n}\) be a a graded algebra and let
\(\rho:sl_2\to\End({\displaystyle \oplus_{n=0}^{2N} A^n})\)
be a given representation compatible with the grading and
the algebra structure (cf.\ Def.\ \ref{Definitioncompatible}).
Such  a representation is given by the image $L$ and $\Lambda$ of
the standard generators of $sl_2$, respectively. We will
consider infinitesimal changes $\rho_\varepsilon$
of $\rho$ such that $L$ deforms to $L+\varepsilon$, where $v\in A^2$ acts
by multiplication.
Prop. \ref{Main1} and \ref{Main2} show that the Lefschetz (or primitive)
decomposition associated to any $sl_2$-representation changes in
a controlled manner when $\rho$ is deformed to such a $\rho_\varepsilon$.

Besides its application to the variation of the space of harmonic forms the precise
formulae also shed some light on the question whether the K\"ahler Lie algebra
of $X$ has the Jordan-Lefschetz property.
The K\"ahler Lie algebra of a compact K\"ahler manifold has recently been
introduced by Looijenga and Lunts \cite{LL}. Instead of looking at just
one K\"ahler class $\omega\in\kk_X$ and its associated $sl_2$-representation
$(L_\omega,\Lambda_\omega,B)$, they study the Lie algebra ${{\mathfrak g}}$
that is generated by all
$L_\omega$ and $\Lambda_\omega$ for $\omega\in\kk_X$. Whereas $L_\omega$ and
$L_{\omega'}$ clearly commute for different $\omega,\omega'\in\kk_X$, the
operators $\Lambda_\omega$ and $\Lambda_{\omega'}$ in general do not.
If they do commute then the Lie algebra ${{\mathfrak g}} $ has the Jordan-Lefschetz property
and Lie algebras of this type have been classified in \cite{LL}.
For tori and hyperk\"ahler manifolds the Jordan-Lefschetz property
can be verified, but for no other class of manifolds
it is known to hold. As it turns out, however, the infinitesimal
Jordan-Lefschetz property always holds.
More precisely, we prove (Thm.\ \ref{LefJorThm}):

\bigskip

{\it Let $X$ be a compact K\"ahler manifold and let $\omega$ be a K\"ahler
class on $X$. If $v$ is  a real $(1,1)$-class, then the contraction operators
$\Lambda_\omega$ and $\Lambda_{\omega+\varepsilon v}$ associated with the K\"ahler
class $\omega$ and its infinitesimal deformation $\omega+\varepsilon v$
commute, i.e.\ $[\Lambda_\omega,\Lambda_{\omega+\varepsilon v}]=0$.}

\bigskip

Of course, in order to decide whether the K\"ahler Lie algebra has the
Jordan-Lefschetz property we need to pass from first order deformations
to higher order deformations, which in general will be obstructed.

As an application of the result about the infinitesimal variation
of the space of harmonic forms we study various sets of special
K\"ahler forms. Usually, a K\"ahler class in $\kk_X$ is lifted to 
the unique Calabi-Yau K\"ahler form which is distinguished by the property
that $\omega^N$ is a scalar multiple of the fixed volume form.
Recently, it has turned out that other forms of the Monge-Amp\`ere equation
are interesting as well. Firstly, in \cite{Leung1} the B-field is chosen
to be a real closed $(1,1)$-form $\beta$ such that $(\omega_0+i\beta)^N$
is a scalar multiple of the volume form. Secondly, Leung \cite{Leung2}
has studied a modified version of the Hermite-Einstein equation
which is related to Gieseker-Maruyama stability (rather than to
slope-stability). It is an equation for the powers of $(\omega_0+tR)$,
where $t$ is a small imaginary parameter and $R$ is the curvature.

So we propose to look at the set $\kk^i$ of
all K\"ahler forms $\omega$ such that
$\omega^{N-i}\omega_0^i$ is a scalar multiple of the volume
form. These are the coefficients of $(\omega_0+t\omega)^N$.
For $i=0$ one obtains the set of Calabi-Yau K\"ahler forms and
for $i=N-1$ the set of $\omega_0$-harmonic K\"ahler forms. For
$0<i\leq N-1$ the lift of a K\"ahler class to such a form is still
unique, but might not always exist. The different sets $\kk^i$
are related to each
other. Each of thes sets is suspected to reflect certain geometric properties
of the variety $X$. E.g.\ we will see that rational curves
in the twistor fibres associated to a Ricci-flat K\"ahler form $\omega_0$
on a K3 surface prevent K\"ahler classes from being represented by
$\omega_0$-harmonic K\"ahler forms (Prop.\ \ref{K3}). Of course,
many aspects of the interplay between the geometry and the metric
structure of Ricci-flat manifolds still remain to be exploited.

The paper is organized as follows.
In Sect.\ \ref{slgeneral} we study deformations of compatible
$sl_2$-representations. Deforming the Lefschetz operator $L$ to $L+\varepsilon
v$ induces a variation of the primitive decomposition $\oplus L^jP^{n-2j}$. This
variation is measured by the natural map
$L^jP^{n-2j}\to\oplus_{i\ne j} L^iP^{n-2i}$, which is described
in Prop.\
\ref{Main1}. The result relies on a formula for the primitive
decomposition of $L^j(\alpha)\in L^jP^{n-2j}$ with respect to
$L+\varepsilon v$ (cf.\ Prop.\ \ref{Main2}).

In Sect.\ \ref{lefschetz} the results are applied to the standard
$sl_2$-representation on cohomology induced by the choice of a K\"ahler class.
In particular, the infinitesimal Jordan-Lefschetz property is proved (Thm.\
\ref{LefJorThm}).

In Sect.\ \ref{harmonicinf} we introduce the map $h(v,_-):\kh_\omega\to
{\rm Im}(d^*)$ that measures the change of the space of $\omega$-harmonic forms
when $\omega$ is deformed to $\omega+\varepsilon v$. Thm.\ \ref{Mainthm} provides
a precise formula for $h(v,\alpha)$ which involves the exterior derivative of
the primitive components of the product $\alpha v$. The
aforementioned criterion that decides
whether $\alpha$ stays harmonic follows from this.

Sect.\ \ref{kaehlercones} is a sequel to \cite{Huy}.
We introduce the series $\kk^0,\kk^1,\ldots,\kk^{N-1}$
of sets of special K\"ahler forms, where $\kk^0$ is the set of Ricci-flat
K\"ahler forms (denoted by $\tilde\kk_X$ in \cite{Huy}) and $\kk^{N-1}$ is the
set of K\"ahler forms that are harmonic with respect to a fixed K\"ahler form
$\omega_0$. Prop.\ \ref{Huyneu} shows that there is a
relation between the harmonicity
of $\alpha^{N-i}\omega_0^i$ for harmonic $(1,1)$-forms $\alpha$ and
the linearity of $\kk^i$. For the set of Ricci-flat K\"ahler forms
$\kk^0$ we complement a result of \cite{Huy} by showing that there is
a relation between the linearity of $\kk^0$ and the harmonicity of
$\alpha^2\omega^{N-2}$ for harmonic $(1,1)$-forms 
(and not of the top exterior power).

\section{Deformations of compatible $sl_2$-representations}\label{slgeneral}

Let \(A={\displaystyle \oplus_{n=0}^{2N} A^n}\) be a finite dimensional
graded $\IC$-algebra.
We assume that $A$ is $\IZ/2\IZ$-commutative with respect to
the induced $\IZ/2\IZ$-grading
$A=A^{even}\oplus A^{odd}$.
A representation $\rho:sl_2(\IC)\to\End_\IC(A)$ is given by the action
of the three generators
$$L=\left(\begin{array}{cc}
0&0\\
1&0\\
\end{array}\right), {\rm ~~}\Lambda=\left(\begin{array}{cc}
0&1\\
0&0\\
\end{array}\right), {\rm ~~and~~}
B=\left(\begin{array}{cc}
1&0\\
0&-1\\
\end{array}\right).$$

\begin{definition}\label{Definitioncompatible}
--- We say that $\rho$ is compatible with the grading
if $L$ and $\Lambda$ are homogenous of degree 2 and -2, resp., and
$B_{|{A^n}}=(N-n)\cdot\id$. The representation $\rho$ is compatible
with the algebra structure of $A$  if $L(\alpha\beta)=L(\alpha)\beta$
for all $\alpha,\beta\in A$. The representation $\rho$ is called compatible
if it is compatible with the grading and the algebra structure.
\end{definition}

Any representation $\rho$ induces a primitive decomposition of $A$. If $\rho$
is compatible with the grading this decomposition
respects the grading, i.e.\ $A^n=\oplus_{j\geq0}L^jP^{n-2j}$,
where $P^m:=\Ker(\Lambda:A^m\to A^{m-2})$  is the space of primitive
elements. Hence, $A=\oplus_{n=0}^{2N}\oplus_{j\geq0}L^jP^{n-2j}$. Recall that
$L^j:A^n\to A^{n+2j}$ is injective for $j\leq N-n$ and that $P^n=0$ for
$n>N$. Also, $\alpha\in A^n$ for $n\leq N$ is primitive if and only if
$L^{N-n+1}(\alpha)=0$.

\begin{example}\label{exex}---
Let $V$ be a vector space of dimension $2N$ with a complex
structure $I$ which respects a scalar product $\langle~~,~~\rangle$ on $V$.
The exterior algebra $\Lambda^*V^*=\oplus_{n=0}^{2N} \Lambda^nV^*$ admits
a natural $sl_2$-representation for which $L$ is given by
multiplication with the K\"ahler form $\omega:=\langle
I(~),~~\rangle\in\Lambda^2V^*$. This is the natural $sl_2$-representation
on hermitian exterior algebras.
\end{example}

Let us consider an infinitesimal deformation $\rho_\varepsilon$ of $\rho$
given by an infinitesimial deformation of the identity $\id\in\End_\IC(A)$,
i.e.\ $\id+\varepsilon\varphi:A[\varepsilon]\to A[\varepsilon]$,
where $\varphi\in \End_\IC(A)$.
More precisely, $\rho_\varepsilon(X)
=(\id+\varepsilon\varphi)^{-1}\cdot\rho(X)\cdot(\id+\varepsilon\varphi)=\rho(X)
+\varepsilon[\rho(X),\varphi]$ for any $X\in sl_2$.
Then the deformation of the generators of $sl_2$ are of the
form $L_\varepsilon:=L+\varepsilon[L,\varphi]$,
$\Lambda_\varepsilon=\Lambda+\varepsilon[\Lambda,\varphi]$,
and $B_\varepsilon=B+\varepsilon[B,\varphi]$.
Of course, as $sl_2$-representations are rigid, the
deformation $\rho_\varepsilon$ as abstract $sl_2$-representation
is isomorphic to the trivial deformation of $\rho$, but the action on the fixed
vector space $A$ in general changes.

If $\rho$ is compatible with the grading, then $\rho_\varepsilon$
is compatible with the grading if and only if $\varphi$ is homogenous
of degree $0$. Then, $B_\varepsilon=B$ and $[L,\varphi]$, respectively
$[\Lambda,\varphi]$ are of degree $2$ and $-2$, respectively.
Deforming $\rho$ this way induces a variation of the decomposition
$A=\oplus_{n=0}^{2N}\oplus_{j\geq0}L^jP^{n-2j}$.
We are interested in two aspects of this.

\bigskip

{\bf i)} The variation of $L^jP^{n-2j}$ as a subspace of $A^n$ is determined by the natural linear map induced by $\varphi$:
$$\tilde\varphi :L^jP^{n-2j}\to\oplus_{i\ne j}L^iP^{n-2i}.$$
{\bf ii)} Any element $L^j(\alpha)\in L^jP^{n-2j}$ admits a primitive
decomposition with respect to $\rho_\varepsilon$.

\bigskip

In general, as $\varphi$ could be any homomorphism of degree $0$, nothing
can be said about {\bf i)} and {\bf ii)}. But as soon as the algebra
structure of $A$ is taken into account the situation  becomes
more interesting.

Let us assume that $\rho$ is a compatible representation on $A$ and that the
deformation $\rho_\varepsilon=\rho+\varepsilon[\rho,\varphi]$ has the property
 that the operator $[L,\varphi]$ of degree two is given by multiplication
with a vector $v\in A^2$, i.e.\ $L_\varepsilon=L+\varepsilon v$.
It is easy to see that such a deformation $\rho_\varepsilon$, where
$\varphi$ is homogenous of degree $0$, is also compatible.
In order to describe the primitive decomposition of $A$ with respect
to $\rho_\varepsilon$ one has in principle to express
$\Lambda_\varepsilon$ in terms of $\Lambda$ and $v$. As there is no easy
way to write this down, we are going to compute the primitive
decomposition directly in terms of $L$ and $v$.

In order to formulate the results we need the following notation.

\begin{definition}---
For $v\in A^2$ and a given representation $\rho$ one defines the linear maps
($i\geq0$):
$$Q_v^i:P^m\to P^{m-2i+2}$$
by the primitive decomposition of the product $v\alpha=\sum_{i\geq0} L^iQ^i_v(\alpha)$.
\end{definition}

In fact, for compatible representations most of these maps are trivial:

\begin{lemma}\label{Griffiths}---
Let $\alpha\in P^m$ and let $\alpha v=\sum_{i\geq0}L^i(\beta_i)$
be the primitive decomposition of $\alpha v$, i.e.\ $\beta_i\in P^{m-2i+2}$.
Then $\beta_i=0$ for $i\geq3$. Equivalently, $Q^i_v=0$ for $i\geq3$.
\end{lemma}

\prf If $\alpha\in P^m$, then $L^{N-m+1}(\alpha)=0$. Hence,
$L^{N-m+1}(\alpha v)=L^{N-m+1}(\alpha)v=0$ and, therefore,
$L^{N-m+1+i}(\beta_i)=0$. Since $L^j$ is injective on $P^{m-2i+2}$ for
$j\leq N-m+2i-2$, one obtains $\beta_i=0$ for $i\geq3$.\qed

\bigskip

The next two propositions show that infinitesimal deformations of compatible
representations of the above form, i.e.\ $L_\varepsilon=L+\varepsilon v$
for some $v\in A^2$, are rather special with regard to {\bf i)}
and {\bf ii)}.

\begin{proposition}\label{Main1}---
Let $\rho$ be a compatible $sl_2$-representation
and let $\rho_\varepsilon$ be a compatible infinitesimal
deformation of $\rho$, such that $L_\varepsilon=L+\varepsilon v$,
where $v\in A^2$ acts by multiplication. The variation of the subspace
$L^jP^{n-2j}\subset A^n$ is measured by the natural map 
$$\tilde\varphi: L^jP^{n-2j}\to\oplus_{i\ne j}L^iP^{n-2i}$$
which is given by
$\tilde\varphi(L^j(\alpha))=(N-n+j+1)L^{j+1}Q^2_v(\alpha)-jL^{j-1}Q_v^0(\alpha)$.
\end{proposition}

The description of $\tilde\varphi$ is an immediate consequence of the following

\begin{proposition}\label{Main2}---
Under the assumption of the previous proposition the primitive decomposition
of $L^j(\alpha)$ for $\alpha\in P^{n-2j}$ with respect to $\rho_\varepsilon$ takes the form
$$\begin{array}{rl}
L^j(\alpha)=&L_\varepsilon^{j+1}(\varepsilon(N-n+j+1)Q^2_v(\alpha))\\
&+L^j_\varepsilon(\alpha-\varepsilon((N-n+2j+1)LQ_v^2(\alpha)+jQ^1_v(\alpha)))\\
&-L_\varepsilon^{j-1}(\varepsilon jQ_v^0(\alpha)).\\
\end{array}$$
\end{proposition}

\prf Using $L_\varepsilon^k=L^k+k\varepsilon v L^{k-1}$ one has
for any $\gamma\in A$ the following equality
$$\begin{array}{rl}L^j(\alpha)=&
L_\varepsilon^{j+1}(\varepsilon\gamma)+
L^j_\varepsilon(\alpha-\varepsilon L(\gamma))\\
&-j\varepsilon(L_\varepsilon^{j-1}Q^0_v(\alpha)+L^j_\varepsilon Q^1_v(\alpha)+L^{j+1}_\varepsilon Q^2_v(\alpha)).\\
\end{array}$$
Since $P^m=0$ for $m>N$, we can assume that $n-2j\leq N$. Furthermore,
$L^{N-n+2j+2}:A^{n-2j-2}\to A^{2N-n+2j+2}$ is bijective. Hence, there is a uniquely
defined $\gamma\in A^{n-2j-2}$ such that $L^{N-n+2j+2}(\gamma)=(N-n+2j+1)L^{N-n+2j}(\alpha v)$. Since $\alpha\in P^{n-2j}$, one has
$$L^{N-n+2j+3}(\gamma)=(N-n+2j+1)L^{N-n+2j+1}(\alpha v)=0.$$
Thus, $\gamma$ is primitive with respect to $\rho$
and, therefore, $\varepsilon\gamma$ is primitive
with respect to $\rho_\varepsilon$. Analogously, since
the $Q_v^i(\alpha)$ are primitive with respect to $\rho$,
the $\varepsilon Q_v^i(\alpha)$ are $\rho_\varepsilon$-primitive.
Moreover,
$(\alpha-\varepsilon L(\gamma))$ is seen to be $\rho_\varepsilon$-primitive
by the following argument:

$$\begin{array}{rcl}
L_\varepsilon^{N-n+2j+1}(\alpha-\varepsilon L(\gamma))&=&
L^{N-n+2j+1}(\alpha)-\varepsilon L^{N-n+2j+2}(\gamma)\\
&&+(N-n+2j+1)\varepsilon L^{N-n+2j}(\alpha v)\\
&=&0.
\end{array}$$

To conclude it suffices to show that $\gamma=(N-n+2j+1) Q_v^2(\alpha)$.
This follows from

$$\begin{array}{rcl}
L^{N-n+2j+2}(\gamma)&=&(N-n+2j+1)L^{N-n+2j}(\alpha v)\\
&=&(N-n+2j+1)L^{N-n+2j}(Q_v^0(\alpha)+LQ_v^1(\alpha)+L^2Q^2_v(\alpha))\\
&=&(N-n+2j+1)L^{N-n+2j+2}Q_v^2(\alpha),\\
\end{array}$$
due to $Q^i_v(\alpha)\in P^{n-2j-2i+2}$, and the injectivity
of $L^{N-n+2j+2}$ on $A^{n-2j-2}$.\qed

\begin{examples}\label{primex}---
The two most interesting special cases of the last proposition are $j=0$ and
$n-2j=0$. Let $\rho$ and $\rho_\varepsilon$  be as before.

{\it i)} If $\alpha\in P^n$, then the primitive decomposition of $\alpha$ with
respect $\rho_\varepsilon$ has the form
$$\alpha=L_\varepsilon((N-n+1)\varepsilon Q^2_v(\alpha))+(\alpha-\varepsilon(N-n+1)LQ_v^2(\alpha)).$$

{\it ii)} If $\alpha=1\in A^0$ and $v\in P^2$, then the primitive decomposition
of $L^j(\alpha)$ with respect to $\rho_\varepsilon$ has the form
$$L^j(1)=L^j(\alpha)=L_\varepsilon^j(\alpha)+L_\varepsilon^{j-1}(-j\varepsilon v).$$

\end{examples}

\begin{corollary}
--- A primitive vector $\alpha\in P^n$ stays primitive with respect
to $L_\varepsilon$ if and only if $L^{N-n}(\alpha v)=0$.
\end{corollary}

\prf The vector $\alpha\in P^n$ is $L_\varepsilon$-primitive if and only
if $Q_v^2(\alpha)=0$ due to the above example.
Since $L^{N-n}(\alpha v)=L^{N-n+2}Q_v^2(\alpha)+L^{N-n+1}Q_v^1(\alpha)
+L^{N-n}Q_v^0(\alpha)=L^{N-n+2}Q_v^2(\alpha)$ and $Q_v^2(\alpha)=0$ if and
only if $L^{N-n+2}Q_v^2(\alpha)=0$, the results follows from this.\qed

\section{Variation of the Lefschetz decomposition}\label{lefschetz}

Everything in the last section immediately applies
to the Lefschetz decomposition of the cohomology of
a K\"ahler manifold. Let $X$ be a compact K\"ahler manifold
and let $[\omega]\in\kk_X\subset H^{1,1}(X,\IR)$ be a K\"ahler class.
The Lefschetz operator $L(\alpha):=[\omega]\wedge\alpha$ is part
of an $sl_2$-representation on $H^*(X)$. The Lefschetz decomposition
is the primitive decomposition induced by this representation
$H^*(X)=\oplus_{j\geq0}L^j H^{n-2j}(X)_{prim}$.

The reason for the existence of an $sl_2$-representation on the
cohomology of a K\"ahler manifold are the K\"ahler identities. The 
natural $sl_2$-representations on the hermitian exterior
algebras $\Lambda^*T^*_x$ for all $x\in X$ (cf.\ Example \ref{exex})
are compatible with passing to cohomology, as the Laplacian
commutes with $L$ and $\Lambda$. Although, two $sl_2$-representations
on $\Lambda^*T^*_x$ associated to two K\"ahler structures commute, this 
does not hold in general on the level of cohomology, since the two
Laplacian are different and $L$ and $\Lambda$ for one K\"ahler structure
do not necessarily commute with the Laplacian of the other. However,
the argument immediately shows that all $sl_2$-representations
on the cohomology of a torus do commute, as the cohomology
is isomorphic to the exterior algebra of the tangent space
at some point.

Changing the K\"ahler class $[\omega]\in\kk_X$ changes this decomposition.
Since an infinitesimal variation of $[\omega]\in\kk_X$ is of the form
$[\omega]+\varepsilon v$ for some $v\in H^{1,1}(X,\IR)_{prim}$, we can apply
Prop.\ \ref{Main2}. In particular, we have

\begin{corollary}---
Let $\alpha\in H^{n-2j}(X)_{prim}$.
Then the Lefschetz decomposition of $L^j(\alpha)$ with respect to
$[\omega]+\varepsilon v$ involves only $L_\varepsilon^{j+1}$,
$L_\varepsilon^j$, and $L_\varepsilon^{j-1}$.\qed
\end{corollary}

Let us next come to the infinitesimal Jordan-Lefschetz property. In \cite{LL}
Looijenga and Lunts introduced the K\"ahler Lie algebra ${{\mathfrak g}}$ of a
compact K\"ahler manifold $X$ as the Lie subalgebra of ${{\mathfrak
gl}}(H^*(X))$ generated by $L_\omega$, $\Lambda_\omega$ for all K\"ahler
classes $\omega$. By \cite[Prop.\ 1.6]{LL} the K\"ahler Lie algebra
${{\mathfrak g}}$ is semi-simple. By Prop.\ 2.1 of the same paper the operators
$\Lambda_\omega$ and $\Lambda_{\omega'}$ commute for all different $\omega$,
$\omega'$ if and only if ${{\mathfrak g}}$ has degree $-2$, $0$, and $2$ only.
In this case $({{\mathfrak g}}, [\Lambda,L])$ is called a Jordan-Lefschetz
pair. Due to \cite[Cor.\ 2.6]{LL} Jordan-Lefschetz pairs can be classified.
Complex tori and compact hyperk\"ahler manifolds satisfy the Jordan-Lefschetz
property, but no other series of higher dimensional compact K\"ahler manifolds
seems to be known. Surprisingly, the infinitesimal Jordan-Lefschetz property
always holds true. This is the content of the next theorem.

\begin{theorem}\label{LefJorThm}---
Let $X$ be a compact K\"ahler manifold and let $\omega\in H^{1,1}(X,\IR)$ be a
K\"ahler class. For any $v\in H^{1,1}(X,\IR)$ the operators $\Lambda_\omega$
and $\Lambda_{\omega+\varepsilon v}$ associated to $\omega$ resp.\ its
infinitesimal deformation $\omega+\varepsilon v$ commute, i.e.\
$[\Lambda_\omega,\Lambda_{\omega+\varepsilon v}]=0$.
\end{theorem}

\prf The assertion is a consequence of Prop.\ \ref{Main2}.
The calculation is straightforward, but lengthy.
We only indicate the main steps.

Let us introduce the notation $\Lambda:=\Lambda_\omega$ and
$\Lambda_\varepsilon:=\Lambda_{\omega+\varepsilon v}$.
In order to prove $[\Lambda_\varepsilon,\Lambda]=0$ it suffices to
show that for any primitive class $\alpha\in H^{n-2j}$ with
$L^j\alpha\ne0$ one has $[\Lambda_\varepsilon,\Lambda](L^j\alpha)=0$.

A simple calculation proves $\Lambda L^k\beta=k(N-m-k+1)L^{k-1}\beta$
for any primitive form $\beta$ of degree $m$. Thus
$\Lambda_\varepsilon\Lambda L^j\alpha=j(N-n+j+1)\Lambda_\varepsilon
L^{j-1}\alpha$. Then, the $(\omega+\varepsilon v)$-primitive
decomposition of $L^{j-1}\alpha$ is given by Prop.\ \ref{Main2}
and the resulting expressions $\Lambda_\varepsilon L_\varepsilon^k$
with $k=j, j-1,j-2$ can be computed as before. Thus, eventually
$\Lambda_\varepsilon\Lambda L^j\alpha$ can be expressed as a linear
combination of $\varepsilon L^{j-1}Q_v^2(\alpha)$,
$\varepsilon L^{j-2}Q_v^1(\alpha)$,
$\varepsilon L^{j-3}Q_v^0(\alpha)$, and $L^{j-2}\alpha$.
Analogously, $\Lambda\Lambda_\varepsilon L^j\alpha$ can be computed
by inserting the formula for the primitive decomposition
of $L^j\alpha$ with respect to $\omega+\varepsilon v$ given
by Prop.\ \ref{Main2}. The resulting formula is
again a linear combination of the same terms as before and,
as it turns out, the linear combination is in both
cases the same. This proves $[\Lambda_\varepsilon,\Lambda]=0$.\qed

\section{Infinitesimal change of harmonic forms}\label{harmonicinf}

In this section we apply the results of Sect.\ \ref{slgeneral}
to $sl_2$-representations on the
space of forms. Although the vector space is no longer finite dimensional
the results are still valid, as the representations are induced by the
standard $sl_2$-representation on hermitian exterior algebras. The
Lefschetz decomposition discussed in the last section
only reflects this on the level of cohomology. The new ingredient
on the level of linear algebra in this context
is the Hodge $*$-operator.

Let $X$ be a compact K\"ahler manifold of dimension $N$.
For any K\"ahler form $\omega$ we denote by $\kh_\omega$ the
space of all forms that are harmonic with respect to $\omega$.
Thus, $\kh_\omega=\oplus\kh_\omega^n=\oplus\kh^{p,q}_\omega$.
Considered as a subspace of the space $\ka(X)_{cl}$ of all
closed forms the space $\kh_\omega$ depends on $\omega$.
We shall compute the first order term of this dependence.
A first order deformation of the K\"ahler form
$\omega$ is of the form $\omega_\varepsilon=\omega+\varepsilon v$,
where $v$ is a closed real $(1,1)$-form.

In order to study the variation of $\kh_\omega$ we have to understand the
deformation of the Laplacian $\Delta$ which involves the Hodge operator
$*$ associated with $\omega$. Let us write the Hodge operator
$*_\varepsilon$ associated to $\omega_\varepsilon=\omega+\varepsilon v$
in the form $*_\varepsilon=*+\varepsilon T_v$, where $T_v$
is a certain linear operator. We shall first describe the deformation
of $\kh_\omega$ to $\kh_{\omega_\varepsilon}$ in terms of $T_v$.
Later we will use the results about the deformation of $sl_2$-representations
to express this purely in terms of $v$.

Recall that the Hodge decomposition of a form $\alpha\in\ka(X)$ has the form
$\alpha=\kh_\omega(\alpha)\oplus dGd^*\alpha\oplus d^* Gd \alpha$, where
$\kh_\omega(\alpha)$ denotes the $\omega$-harmonic part of $\alpha$
and $G$ is the Green operator with respect to $\omega$, which
commutes with $d$ and $d^*$. Analogously, one has $\alpha=\kh_{\omega_\varepsilon}
(\alpha) \oplus dG_\varepsilon d^{*_\varepsilon}\alpha\oplus
d^{*_\varepsilon} G_\varepsilon d\alpha$. The infinitesimal variation
of $\kh_\omega\subset\ka(X)_{cl}$ induced by $\omega_\varepsilon=\omega+\varepsilon v$ is determined by the canonical homomorphism
$$\tilde h(v,_-):\kh_\omega\to\ka(X)_{cl}/\kh_\omega={\rm Im}(d).$$
The space of harmonic forms $\kh_\omega$ does not change
when passing from $\omega$ to $\omega_\varepsilon$ if and only if this
map is trivial. More specifically, one has

\bigskip

{\it An $\omega$-harmonic form $\alpha\in\kh_\omega$ is harmonic
with respect to $\omega_\varepsilon$, i.e.\
$\alpha\in\kh_{\omega_\varepsilon}$, if and only if $\tilde h(v,\alpha)=0$.}

\bigskip

The map $\tilde h(v,_-)$ can explicitely be described in terms
of $T_v$ due to the following

\begin{lemma}\label{harmoniceps}
--- Let $\alpha\in\kh_\omega$. Then the $\omega_\varepsilon$-harmonic
part of $\alpha$ is given by $\kh_{\omega_\varepsilon}(\alpha)=\alpha
+\varepsilon dG*dT_v(\alpha)$.
\end{lemma}

\prf
Since $\varepsilon dG \ast dT_v(\alpha) $ is exact, it suffices to
show that $\beta := \alpha + \varepsilon dG\ast
dT_v (\alpha)$ is harmonic with respect to $\omega_\varepsilon$.
Obviously, $\beta$ ist $d$-closed. In order to show that
$d^{\ast_\varepsilon} \beta =0$, one first notes that
$\varepsilon G = \varepsilon G_\varepsilon$.
Also recall that $\Delta G = G\Delta =
{\rm id}$ and $\Delta_\varepsilon G_\varepsilon = G_\varepsilon\Delta_\varepsilon ={\rm id}$ on $d$-exact forms and that $d^*=-*d*$.
\\
Therefore,
$$
\begin{array}{l}
 \Delta_\varepsilon (\alpha+ \varepsilon dG\ast dT_v(\alpha ))=
\Delta_\varepsilon(\alpha ) +\varepsilon \Delta_\varepsilon G_\varepsilon
d \ast dT_v \alpha \\
=dd^{\ast_\varepsilon} \alpha + \varepsilon d \ast dT_v (\alpha)\\
=-d(\ast +\varepsilon T_v) d
(\ast +\varepsilon T_v)(\alpha ) + \varepsilon d \ast dT_v (\alpha )=0.\\
\end{array}
$$
\phantom{i}\hfill$\square$

\bigskip

Thus, the Hodge decomposition of $\alpha$ with respect to
$\omega_\varepsilon$ takes the form
$\alpha=(\alpha+\varepsilon d G*dT_v(\alpha))\oplus
(-\varepsilon dG*dT_v(\alpha))$.

\begin{corollary}\label{descriptionh}---
The map $\tilde h (v,_- ) : \kh_\omega  \lra {\rm Im} (d)$ maps $\alpha $ to
$- d G\ast d T_v (\alpha)$\hfill$\square$
\end{corollary}

\begin{remark}---
The Green operator commutes with $d$ and defines an automorphism of
${\rm Im} (d)$. Moreover, $d: {\rm Im} (d^\ast )\to {\rm Im} (d)$ is
bijective. So, without loosing any information we may replace $\tilde h (v,_- )$
by the map $h(v,_- ): \kh_\omega \to {\rm Im}(d^\ast )$ that maps $\alpha $
to $\ast dT_v(\alpha) $. Note that $h(v,_- )$ is a map of degree $(-1)$.
In fact, we could also consider $dT_v(\alpha) $ without loosing
information, but the degree of this map depends on $\alpha $ and on $N$.
\end{remark}

Next we shall compute $T_v(\alpha) $ in terms of the product
$\alpha v$. Throughout, we will make use of the
following formula due to A.\ Weil (cf. \cite{Weil, Wells}):
\\
Let $\alpha $ be a form of type $(p,q)$. Assume that $\alpha $ is
$\omega$-primitive and that $(p+q)+r \le N$. Then
\\
$$\ast (\omega^r \alpha ) = (-1)^{{\frac{(p+q)(p+q+1)}{2}}} \cdot {\frac{r!\
i^{p-q}}{{(N - (p+q)-r)!}}} \cdot \omega^{N-(p+q)-r} \alpha \ .
\eqno{(1)}$$

In fact, if both sides are correctly interpreted the formula also
holds for $(p+q)+r>N$. Namely, $\omega^r \alpha = 0$, since $\alpha $
is primitive and $r> N-(p+q)$, and $\omega^{N-(p+q)-r} \alpha = 0$, as
the exponent is negative.

The idea to compute $T_v(\alpha)$ is the following. Using Prop.\ \ref{Main2} we
can compute the primitive decomposition of any $\alpha\in\kh_\omega$ with
respect to $\omega_\varepsilon$. Formula (1) allows one
to compute $*_\varepsilon$
by applying it to each of the summands in the
$\omega_\varepsilon$-primitive
decomposition of $\alpha$.
On the other hand, $*_\varepsilon(\alpha)=*(\alpha)+\varepsilon T_v(\alpha)$
can also be computed by applying (1) to $*$. This gives a formula for
$T_v(\alpha)$.

\begin{example}--- Let $\alpha$ be an $\omega$-primitive form
of type $(1,1)$. By \ref{primex} the $\omega_\varepsilon$-primitive
decomposition of $\alpha$ is given by
$\alpha=\omega_\varepsilon((N-1)\varepsilon Q_v^2(\alpha))
+(\alpha-\varepsilon(N-1)\omega Q_v^2(\alpha))$. The function $f=Q_v^2(\alpha)$
is determined by $Q_v^2(\alpha)\omega^N= \alpha v\omega^{N-2}$. Thus,
$f=(\alpha v\omega^{N-2})/\omega^N$. Hence, the $\omega_\varepsilon$-primitive
decomposition of $\alpha$ is
$$\alpha=\omega_\varepsilon((N-1)\varepsilon\left(\frac{\alpha
v\omega^{N-2}}{\omega^N}\right))+(\alpha-\varepsilon(N-1)\omega\left(\frac{\alpha
v\omega^{N-2}}{\omega^N}\right)).$$ By (1):

$$\begin{array}{rcl}
*_\varepsilon\alpha&=&*_\varepsilon(\omega_\varepsilon((N-1)\varepsilon\left
(\frac{\alpha v\omega^{N-2}}{\omega^N}\right)))+*_\varepsilon
(\alpha-\varepsilon(N-1)\omega\left(\frac{\alpha v\omega^{N-2}}{\omega^N}\right))\\
&=&\frac{N-1}{(N-1)!}\varepsilon\omega_\varepsilon^{N-1}\left(\frac{\alpha v\omega^{N-2}}{\omega^N}\right)-\frac{1}{(N-2)!}\omega_\varepsilon^{N-2}(\alpha-\varepsilon(N-1)\omega\left(\frac{\alpha v\omega^{N-2}}{\omega^N}\right)).\\
\end{array}$$

On the other hand, $*_\varepsilon\alpha=*\alpha+\varepsilon T_v(\alpha)=-\frac{1}{(N-2)!}\omega^{N-2}\alpha+\varepsilon T_v(\alpha)$.
Hence,
$$ T_v(\alpha)=-\frac{1}{(N-3)!} \omega^{N-3} \alpha v +
\frac{N}{(N-2)!} \left( {\frac{\alpha v\omega^{N-2}}{\omega^N}}\right)
\omega^{N-1} .$$

- For the second example let us assume that $v$ is primitive. Then
for the $(1,1)$-form $\alpha = \omega$ the primitive
decomposition with respect to $\omega_\varepsilon = \omega + \varepsilon v$
takes the form $\omega = (\omega +\varepsilon v)-\varepsilon v$.
Indeed, $\varepsilon v$ is $\omega_\varepsilon$-primitive, for
$\varepsilon v(\omega + \varepsilon v)^{N-1} = \varepsilon v
\omega^{N-1}=0$.
\\
Thus,
$$
\ast_\varepsilon\omega = \ast_\varepsilon \omega_\varepsilon-
\varepsilon \ast_\varepsilon v
= \frac{1}{(N-1)!} {\omega_\varepsilon}^{N-1} +
{\frac{\varepsilon}{(N-2)!} }{\omega_\varepsilon}^{N-2}v \ .
$$

Using, $\ast_\varepsilon = \ast + \varepsilon T_v$ this yields

$$T_v( \omega) = {\frac{(N-1)}{(N-1)!} }\omega^{N-2} v + \frac{1}
{(N-2)!} \omega^{N-2} v =\frac{2}{(N-2)!} \omega^{N-2}v\ .$$
\end{example}

\bigskip

The general formula is provided by the following

\begin{proposition}\label{Tvdescription} ---
Let $\alpha$ be $\omega$-primitive of type $(p-j,q-j)$ and let
$\alpha v=\beta_0 +\omega \beta_1+\omega^2\beta_2$ be the $\omega$-primitive
decomposition of $\alpha v$. Then
$$T_v (\alpha\omega^j)= c'_0 \beta_0 \omega^{N-n+j-1} +
c_1 \beta_1 \omega^{N-n+j} +c_2 \beta_2 \omega^{N-n+j+1} +
c_0 \alpha  v \omega^{N-n+j-1}\ ,$$
\noindent
where $c_0 = -\eta \cdot (N-n+j)$, $c_1 = \eta \cdot j$,
$c_2 = \eta \cdot (N-n+3j+2)$, $\eta :=
(-1)^{ {\frac{(n-2j-2)(n-2j-1)}{2}}} \  {\frac{j! i^{p-q}}{(N-n+j)!}}$,
$n=p+q$, and
$c'_0 =c_0$ for $j>0$ and $c_0'=0$ for $j=0$.
\end{proposition}

\prf Using $(\ast +\varepsilon T_v)(\alpha \omega^j ) = \ast_\varepsilon
(\alpha \omega^j )$ and the primitive decomposition of $\alpha \omega^j$ with
respect to $\omega_\varepsilon$ given by Prop.\ \ref{Main2}, the claim is
proven by applying (1). By Prop.\ \ref{Main2} we have the
$\omega_\varepsilon$-primitive decomposition
$$\begin{array}{rl}
\omega^j\alpha=&(N-n+j+1)\omega_\varepsilon^{j+1}(\varepsilon\beta_2)\\
&+\omega_\varepsilon^j(\alpha-\varepsilon((N-n+2j+1)\omega\beta_2+j\beta_1))\\&-j\omega_\varepsilon^{j-1}(\varepsilon\beta_0).\\
\end{array}$$
Applying $*_\varepsilon$ to these three terms gives:
$$\begin{array}{rcl}
\ast_\varepsilon {\omega_\varepsilon}^{j-1} (\varepsilon (-j \beta_0)) &=&
c_0 \varepsilon \beta_0 {\omega_\varepsilon}^{N-n+j-1} \quad \hbox{ for } j>0 ,\\
\ast_\varepsilon {\omega_\varepsilon}^j (\alpha - \varepsilon ((N-n+2j+1)\beta_2 \omega +j\beta_1))
&=&-\eta {\omega_\varepsilon}^{N-n+j}(\alpha-\varepsilon ((N-n+2j+1)\beta_2\omega +
j\beta_1 ))\\
{\rm and~~~~~}
\ast_\varepsilon {\omega_\varepsilon}^{j+1} (\varepsilon (N-n+j+1)\beta_2 ) &=&
\eta (j+1) {\omega_\varepsilon}^{N-n+j+1} (\varepsilon \beta_2 ) \ .
\end{array}$$
The scalars $c_1 ,c_2$ can be easily computed from this.\hfill$\square$

\bigskip

\begin{remark}---
Note that if $N-n+j < 0$, then $\alpha \omega^j =0$. In this case
the right hand side is also interpreted as zero. For
$N-n+j =0$ the right hand side reduces to $c_1 \beta_1$.

- For $j=0$ the formula yields
$$T_v (\alpha )= (-1)^{\frac{n(n+1)}{ 2}} \ {\frac{i^{p-q}}{{(N-n)!}}}
\left( (N-n)\alpha v \omega^{N-n-1} - (N-n+2)\beta_2
\omega^{N-n+1}\right). $$
- For $p=q=1$ we get back the formula in the example, since $\beta_2=(\alpha v\omega^{N-2})/\omega^N$.

- For $\alpha =1$ and $v$ primitive we find
$T_v (\omega^j )= {\frac{2j!}{(N-j-1)!}}v\omega^{N-j-1}$.
\end{remark}

Passing from $T_v(\alpha \omega^j)$ to $h(v,\alpha \omega^j)=*dT_v(\alpha\omega^j)$ actually simplifies the formula due to the following

\begin{lemma}\label{beta0harmonic}
--- Let $\alpha$ be a closed primitive form and
$\alpha v=\beta_0+\beta_1\omega+\beta_2\omega^2$ be the
primitive decomposition of the product $\alpha v$, where
$v$ is a closed form of degree two. Then the primitive decomposition
 of $d\beta_0$, $d\beta_1$, and $d\beta_2$ are of the form:
$d\beta_1=\delta_0+\delta_1\omega$, $d\beta_0=-\delta_0\omega$, and
$d\beta_2=-\delta_1$. In particular, $d\beta_1=0$ if and only if
$d\beta_0=d\beta_2=0$.
\end{lemma}

\prf 
Let $d\beta_0=\sum\delta_k^0\omega^k$,
$d\beta_1=\sum\delta_k\omega^k$, and 
$d\beta_2=\sum\delta_k^2\omega^k$ be the primitive decomposition
of $d\beta_0$, $d\beta_1$, and $d\beta_2$, respectively.
We use the following three equations to deduce the result:
{\it i)} $d(\alpha v)=0$,
{\it ii)} $(d\beta_1+d\beta_2\omega)\omega^{N-\deg(\alpha)}=0$, 
and {\it iii)} $d\beta_2\omega^{N-\deg(\alpha)+2}=0$.

Since $\alpha$ and $v$ are closed one has {\it i)}, which
yields 
$$\sum\delta_k^0\omega^k+\sum\delta_k\omega^{k+1}+
\sum\delta_k^2\omega^{k+2}=0.$$
Hence, $\delta_0^0=0$, $\delta_1^0=-\delta_0$, and
$\delta_{k+2}^0+\delta_{k+1}+\delta_k^2=0$ for all $k\geq0$. 
If $\deg(\alpha)\geq N-1$ then $\beta_0=0$. If $\deg(\alpha)\leq N-1$, we use 
that $\beta_0$ is of degree $\deg(\alpha)+2$ to conclude that
$\beta_0\omega^{N-\deg(\alpha)-1}=0$. This yields {\it ii)} and therefore
$$\sum\delta_k\omega^{N-\deg(\alpha)+k}+\sum\delta_k^2\omega^{N-\deg(\alpha)+k+1}=0.$$
Since $\delta_{k+1}+\delta_k^2$ is primitive of degree
$\deg(\alpha)-1-2k$ we obtain $\delta_{k+1}+\delta_k^2=0$ for all
$k\geq0$. This already proves $d\beta_0=\delta_1^0\omega=-\delta_0\omega$.

Using that $\beta_1$ is primitive of degree $\deg(\alpha)$ equation
{\it ii)} gives {\it iii)}, i.e.\
$\sum\delta^2_k\omega^{N-\deg(\alpha)+2+k}=0$.
Since $\delta_k^2$ is primitive of degree $\deg(\alpha)-1-2k$, one obtains
$\delta_k^2=0$ for $k>0$.\qed

\begin{lemma}\label{alphavharmonicbeta1}
--- Let $\alpha$ be a closed primitive form of pure type and
let $v$ be closed of type (1,1). Then $\alpha v$ is harmonic if
and only if $\beta_1$ is closed. As before,
$\beta_1$ is the primitive form defined by the primitive decomposition
$\alpha v=\beta_0+\beta_1\omega+\beta_2\omega^2$.
\end{lemma}

\prf Since $X$ is K\"ahler, $\alpha v$ is harmonic if and only if $\beta_0$,
$\beta_1$, and $\beta_2$ are harmonic. On the other hand,
if $d\beta_1=0$, then also $d\beta_0=d\beta_2=0$. Since the forms
$\beta_i$ are primitive of pure type, this implies
$d*\beta_i=0$, i.e.\ they are harmonic.
(In general, a closed primitive form of pure type is harmonic. In particular,
$\alpha$ is harmonic from the very beginning.)\qed

This leads to the final formula.

\begin{theorem}\label{Mainthm}
--- Let $\alpha$ be $\omega$-primitive and $\omega$-harmonic
of type $(p-j,q-j)$ with $\alpha\omega^j\ne0$. If $\alpha v=\beta_0+\beta_1
\omega+\beta_2\omega^2$ is the primitive decomposition of $\alpha v$
and $d\beta_1=\delta_0+\delta_1\omega$ is the primitive decomposition of
$d\beta_1$, then
$$h(v,\alpha\omega^j)=\lambda_1*(\delta_0\omega^{N-n+j})+
\lambda_2*(\delta_1\omega^{N-n+j+1}),$$
where $\lambda_1=\eta(N-n+2j)$ for $j>0$ and $\lambda_1=0$
for $j=0$, $\lambda_2=-\eta(N-n+2j+2)$, $\eta$ as before, and
$n=p+q$. 

\end{theorem}

\prf By definition $h(v,\alpha\omega^j)=*dT_v(\alpha\omega^j)$ and
by Prop.\ \ref{Tvdescription}
$$\begin{array}{rcl}
*dT_v(\alpha\omega^j)&=&*(c_0' d\beta_0\omega^{N-n+j-1}+
c_1d\beta_1\omega^{N-n+j}+c_2d\beta_2\omega^{N-n+j+1})\\
&=&(c_1-c_0')*(\delta_0\omega^{N-n+j})+(c_1-c_2)*(\delta_1\omega^{N-n+j+1}).\\
\end{array}$$
The calculation of $\lambda_1:=c_1-c_0'$ and $\lambda_2:=c_1-c_2$ is
straightforward.
\qed

\begin{corollary}\label{harmonicav}
--- If $j=0$ the form $\alpha$ is harmonic
with respect to $\omega_\varepsilon=\omega+\varepsilon v$
if and only if $\alpha v\omega^{N-n}$ is $\omega$-harmonic.
If $j>0$  the form $\alpha\omega^j$ is harmonic
with respect to $\omega_\varepsilon=\omega+\varepsilon v$
if and only $\alpha v$ is $\omega$-harmonic.
\end{corollary}

\prf Let first $j=0$. Since $\lambda_1=0$ for $j=0$ one has
$h(v,\alpha)=0$ if and only if $\delta_1\omega^{N-n+1}=0$.
Since $\delta_1$ is of degree $n-1$ the latter is equivalent
to the vanishing of $d\beta_2$. Now, $\beta_2$ is closed if and only
if $\beta_2$ is harmonic if and only if $\beta_2\omega^{N-n+2}$
is harmonic, but the latter is just $\alpha v\omega^{N-n}$.

In general, if $\alpha v$ is harmonic, then $\delta_0=\delta_1=0$ by lemma 
\ref{alphavharmonicbeta1} and hence $h(v,\alpha\omega^j)=0$.
Conversely, if $h(v,\alpha\omega^j)=0$, then
$\delta_1=0$, because $\delta_1$ is primitive of degree $n-2j-1$.
For $j>0$ the same argument shows $\delta_0=0$.\qed

\section{Non-linear K\"ahler cones}\label{kaehlercones}

Let $X$ be a compact K\"ahler manifold of dimension $N$ with a fixed K\"ahler
form $\omega_0$. By $\kk_X\subset H^{1,1}(X,\IR)$ we denote the K\"ahler
cone, i.e.\ the set of all K\"ahler classes. Then $[\omega_0]\in \kk_X$.
Clearly, $\kk_X$ is an open convex cone. By definition, every class
in $\kk_X$ can be represented by a K\"ahler form, but a priori there
is no canonical choice.

\begin{definition}--- For $i=0,\ldots, N-1$ we denote by $\kk^i=\kk^i_{(X,\omega_0)}$ the connected component containing $\omega_0$ of the set of all
K\"ahler forms $\omega$, such that $\omega^{N-i}\omega_0^i=c\cdot\omega_0^N$
for some scalar constant $c$.
\end{definition}

\begin{remark}--- The set $\kk^0=\kk^0_{(X,\omega_0)}$ only depends on the
volume form $\omega_0^N$ and not on $\omega_0$ itself. In \cite{Huy}
this set was denoted by $\tilde\kk_X$. Due to results of Calabi and
Yau the natural projection $\kk^0\to\kk_X$, $\omega\mapsto[\omega]$
is bijective. If $\omega_0$ is Ricci-flat, then $\kk^0$ is the set of all
Ricci-flat K\"ahler forms.

- In general, the canonical projection $\kk^i\to\kk_X$ is injective.
This is proved using the original  argument of Calabi \cite{Calabi}:
If $\omega,\omega'\in\kk^i$ with $[\omega]=[\omega']$,
then
$\omega-\omega'=dd^c\varphi$ and $dd^c\varphi(\omega^{N-i-1}+
\omega^{N-i-2}\omega'+\ldots+{\omega'}^{N-i-1})\omega_0^i=0$. The positivity of $\omega$,
$\omega'$, and $\omega_0$ yields $\varphi=0$.

- In general, the projection $\kk^i\to\kk_X$ need not be surjective for
$i>0$ (see the example below). Presumably, taking the
connected component in the definition of $\kk^i$ is superflous.
It certainly is for $i=0$ and $i=N-1$.

- The only `linear' cone is $\kk^{N-1}$. Indeed,
$\kk^{N-1}$ consists of all K\"ahler forms $\omega$ such that
$\omega\omega_0^{N-1}=c\cdot\omega_0^{N}$. Thus, $\omega-c\omega_0$ is
closed and $\omega_0$-primitive, hence $\omega_0$-harmonic. Therefore,
$\kk^{N-1}$ is an open subset of $(\kh^{1,1}_{\omega_0})_\IR$.
In some imprecise sense the sequence $\kk^0,\kk^1,\ldots,\kk^{N-1}$ is
a stepwise linearization of the Calabi-Yau cone $\kk^0$.

\end{remark}

\begin{example}---
Often, for a class $[\alpha]\in\partial\kk_X$ there exists a curve $C\subset X$
such that $\int_C\alpha=0$. Assume that $\kk^{N-1}\to\kk_X$ is surjective
and let $\alpha\in\partial\kk^{N-1}\subset\kh^{1,1}_{\omega_0}$ be
representing $[\alpha]$. Since $\alpha$ is in the boundary of $\kk^{N-1}$
it is semi-positive definite. On the other hand, $\int_C\alpha=0$.
Therefore, $\alpha_{|_C}=0$. Conversely, if we want to find examples
for which $\kk^{N-1}\to\kk^0$ is not surjective, then we have
to look for $[\alpha]$ such that $\alpha_{|_C}\ne0$.

Let $C$ be a curve of genus two and $\varphi_1,\varphi_2\in H^{1,0}(C)$ be an
orthogonal base, i.e.\ $\int\varphi_1\bar{\varphi}_2=0$. Consider the
form
$\alpha=\varphi_1\times\bar{\varphi}_2+\bar{\varphi}_1\times\varphi_2$ on
$X:=C\times C$, which is harmonic with respect to $\omega\times\omega$
for any K\"ahler
form $\omega$ on $C$. Then $\alpha$ is trivial on the fibres of the two
projections and $\int_\Delta\alpha=0$, where $\Delta\subset C\times C$ is the
diagonal. On the other hand, we may assume that
$\alpha_{|_\Delta}=\varphi_1\wedge\bar\varphi_2+
\bar\varphi_1\wedge\varphi_2$ is not trivial. Indeed,
we may change $\varphi_1$ by a complex scalar $\lambda$ and
if for all $\lambda$ one has
$0=\lambda\varphi_1\wedge\bar\varphi_2+\bar\lambda
\bar\varphi_1\wedge\varphi_2$ then $\varphi_1\wedge\bar\varphi_2=0$,
but $\varphi_1$ and $\varphi_2$ vanish at different points.
Therefore, $\alpha$  is not semi-positive definite. Since obviously
$[\alpha]\in\partial\kk_X$, the form
$\alpha$ itself cannot be in $\partial\kk^{N-1}$. Hence, $\kk^{N-1}\to\kk_X$
is not surjective.

\end{example}

Here is another geometric relevant example which shows that 
the existence of certain subvarities of the (deformed) manifold
influences the size of the set of all 
$\omega_0$-harmonic K\"ahler forms.

\begin{proposition}\label{K3}
--- Let $\omega_I$ be a Ricci-flat K\"ahler form on a K3 surface $(X,I)$,
where $I$ is a fixed complex structure on $X$, and let
$\IP:=\{aI+bJ+cK| a^2+b^2+c^2=1\}$ be the induced twistor family
of complex structures.
Let $\alpha\in \kk_X$ be a K\"ahler class on $(X,I)$ and
$\omega_I+\alpha'$ be its primitive decomposition with
respect to $\omega_I$. If for some complex
structure $\lambda\in\IP$ there exists a smooth rational curve $C$ in
$(X,\lambda)$ such that $C.(\omega_\lambda+\alpha')<0$, then
$\alpha$ cannot be represented by an $\omega_I$-harmonic
K\"ahler form.
\end{proposition}

\prf We denote by $\alpha$ and $\alpha'$ also the harmonic
representatives of $\alpha$ and $\alpha'$, respectively.
Then $\alpha'$ is of pure type $(1,1)$ for any complex structure
of the form $\lambda=aI+bJ+cK$. Let us first show that
$(\omega_\lambda+\alpha')^2$ is independent
of $\lambda$ at any point of $X$.
Indeed, $(\omega_\lambda+\alpha')^2=\omega_\lambda^2+{\alpha'}^2=
\omega_I^2+{\alpha'}^2=\alpha^2$, as $\alpha'$ is
$\omega_\lambda$-primitive for any $\lambda$.
Therefore, if the harmonic representative of $\alpha$ is in fact positive
definite, then $\omega_\lambda+\alpha'$ is a positive
definite form on $(X,\lambda)$ for all $\lambda\in\IP$, as otherwise the
square would have to vanish for some triple $\lambda$ at at least
one point of $X$. But then a rational curve as above cannot exist.
Thus, $\alpha$ cannot be represented by an $\omega_I$-harmonic K\"ahler form.
\qed 

\bigskip

In other words, the part of the K\"ahler cone $\kk_X$ of the
K3 surface $(X,I)$ that can be represented by harmonic
K\"ahler forms stays constant in the twistor family.
In particular, this excludes the existence of a series
of smooth rational curves $C_i$ on twistor fibres $(X,\lambda_i)$ with
$\int_{C_i}\omega_{\lambda_i}\to0$. It is an open and interesting question whether
any K\"ahler class that stays K\"ahler in the twistor family can actually
be represented by an harmonic K\"ahler form.

\bigskip

Let us come back to the general situation.
The description of the tangent space of $\kk^0=\tilde\kk_X$ in
\cite{Huy} can easily be adapted to describe the tangent space of $\kk^i$
at the point $\omega_0\in\kk^i$.

\begin{lemma}---
At the point $\omega_0\in\kk^i$ the tangent space is
$T_{\omega_0}\kk^i=(\kh^{1,1}_{\omega_0})_\IR$, i.e.\ any first order
deformation of $\omega_0$ within $\kk^i$
is given as $\omega_0+\varepsilon v$ with $v\in(\kh^{1,1}_{\omega_0})_\IR$.
\end{lemma}

\prf Indeed, $(\omega_0+\varepsilon v)^{N-i}\omega_0^i=
\omega_0^N+(N-i)\varepsilon\omega_0^{N-1}v$. Hence,
for some scalar $\lambda$ the form
$v-\lambda\omega_0$ is $\omega_0$-primitive and closed. Thus, $v$ is
$\omega_0$-harmonic,\qed

\bigskip

Thus, the sets $\kk^i$ have contact at $\omega_0$ of
order two. The tangent space
of $\kk^0$ can be described at every point $\omega\in\kk^0$ as the space
of all real $\omega$-harmonic $(1,1)$-forms.
The description of $T_\omega\kk^i$ for $0<i<N-1$ and $\omega\ne\omega_0$
is less clear, it mixes $\omega$-harmonicity and  $\omega_0$-harmonicity.

We next want to know under what circumstances the `non-linear'
cone $\kk^i$ is linear, i.e.\ when it is contained in its tangent
space $T_{\omega_0}\kk^i=(\kh^{1,1}_{\omega_0})_\IR$. The following
proposition generalizes \cite[Prop.\ 2.3]{Huy}

\begin{proposition}\label{Huyneu}---
The cone $\kk^i$ is linear if and only if $\alpha^{N-i}\omega_0^i$ is
$\omega_0$-harmonic for all $\alpha\in\kh^{1,1}_{\omega_0}$.
\end{proposition}

\prf If $\kk^i$ is linear, then for all $\alpha$ in the open
set $\kk^i\subset(\kh^{1,1}_{\omega_0})_\IR$ the form
$\alpha^{N-i}\omega_0^i$ is $\omega_0$-harmonic. But then this
holds for all $\alpha\in \kh^{1,1}_{\omega_0}$. Conversely, if
$\alpha^{N-i}\omega_0^i$ is $\omega_0$-harmonic for all $\alpha\in\kh^{1,1}_{\omega_0}$, then $\kk^i$ intersects $(\kh^{1,1}_{\omega_0})_\IR$ in an open
subset. Hence, $\kk^i\subset(\kh^{1,1}_{\omega_0})_\IR$.\qed

\bigskip

\begin{corollary}---
If $\kk^i$ is linear, then $\kk^{i+1}$ is linear.
\end{corollary}

\prf If $\kk^i$ is linear then $\alpha^{N-i}\omega_0^i$ is
$\omega_0$-harmonic for all $\alpha\in\kh^{1,1}_{\omega_0}$. Hence,
$(\alpha+t\omega_0)^{N-i}\omega_0^i$ is $\omega_0$-harmonic for all $t$.
Then, also the linear coefficient of this polynomial in $t$,
which is $\alpha^{N-i-1}\omega_0^{i+1}$,
is $\omega_0$-harmonic for all $\alpha\in\kh^{1,1}_{\omega_0}$. Thus,
$\kk^{i+1}$ is linear.\qed

\bigskip

In particular, the linearity of $\kk^0$ implies the linearity of all the
other sets $\kk^i$. In this case, $\kk^0=\kk^1=\ldots=\kk^{N-1}$.
In this sense, the sequence $\kk^0$, $\kk^1$, \ldots, $\kk^{N-1}$
goes from the most curved set $\kk^0$ to the linear set $\kk^{N-1}$.

The following result is another criterion for the linearity
of $\kk^0$.

\begin{proposition}---
The set $\kk^0$ is linear, i.e.\ $\kk^0=\kk^1=\ldots=\kk^{N-1}$,
if and only if $\kh^{1,1}_\omega=\kh^{1,1}_{\omega_0}$ for all
$\omega$ in a neighbourhood of $\omega_0$ in
$\kk^{N-1}\subset(\kh^{1,1}_{\omega_0})_\IR$.
\end{proposition}

\prf If $\kk^0$ is linear, then $\kh^{1,1}_\omega=\kh^{1,1}_{\omega_0}$
for all $\omega \in\kk^0=\kk^{N-1}$. Conversely, if $
\kh^{1,1}_\omega=\kh^{1,1}_{\omega_0}$ for
all $\omega\in(\kh^{1,1}_{\omega_0})_\IR$ close to
$\omega_0$, then $\kh^{N,N}_{\omega+\varepsilon v}=\kh^{N,N}_\omega$
for all $v\in T_\omega\kk^{N-1}=(\kh^{1,1}_\omega)_\IR$.
Indeed, $\kh^{N,N}_{\omega_+\varepsilon v}=\IC(\omega+\varepsilon v)^N=
\IC(\omega^N+(N-1)\varepsilon\omega^{N-1}v)$ and $\omega^{N-1}v$ is
$\omega$-harmonic. Thus, the map $\kk^{N-1}\to\IP(\ka^{N,N}(X))$, $\omega
\mapsto\kh^{N,N}_\omega$ has vanishing differential at every point in
an open neighbourhood of $\omega_0$ and is, therefore, constant.
In particular, for
an open subset of $(\kh^{1,1}_{\omega_0})_\IR$ the top exterior power
$\alpha^N$ is $\omega_0$-harmonic. Therefore, the intersection
of $\kk^0$ and $(\kh^{1,1}_{\omega_0})_\IR$ contains an open
subset. This implies $\kk^0\subset(\kh^{1,1}_{\omega_0})_\IR$,
i.e.\ $\kk^0$ is linear.\qed

\bigskip

The results of Sect.\ \ref{harmonicinf}
yield

\begin{proposition}---
The following conditions are equivalent.

{\it i)} The set $\kk^0$ is linear.

{\it ii)} For all $\alpha\in\kh^{1,1}_{\omega_0}$ the form $\alpha^N$ is
$\omega_0$-harmonic.

{\it iii)} For all $\omega\in\kk^0$ and all $\alpha\in\kh^{1,1}_\omega$
the form $\alpha^N$ is $\omega$-harmonic.

{\it iv)} For all $\omega\in\kk^0$ and all $\alpha\in\kh^{1,1}_\omega$ the
form $\alpha^2\omega^{N-2}$ is $\omega$-harmonic.

\end{proposition}

\prf The equivalence of {\it i)} and {\it ii)}
was shown in \cite{Huy}. Clearly, {\it iii)} implies
{\it ii)}. On the other hand, if $\kk^0$ is linear, then $\kh^{1,1}_\omega=\kh^{1,1}_{\omega_0}$ and $\kh^{N,N}_\omega=\kh^{N,N}_{\omega_0}$ for all
$\omega\in\kk^0$. Hence, {\it i)} implies {\it iii)}. By Corollary \ref{harmonicav}
the form
$\alpha^2\omega^{N-2}$ is $\omega$-harmonic for all $\alpha\in \kh^{1,1}_\omega$ if
and only if $\kh^{1,1}_\omega=\kh^{1,1}_{\omega+\varepsilon v}$ for
all $v\in(\kh^{1,1}_\omega)_\IR$. Hence, {\it iv)} holds if and only if the
Gauss map of the embedding $\kk^0\subset\ka^{1,1}(X)_{cl}$ has everywhere
vanishing differential, i.e.\ $T_\omega\kk^0$ is constant.\qed

\bigskip

Note that condition {\it iv)} for $\omega=\omega_0$ is equivalent to
$\kk^{N-2}$ being linear (Prop.\ \ref{Huyneu}). Thus, $\kk^0$ is linear
if and only if $\kk^{N-2}_{(X,\omega)}$ is linear for any $\omega\in\kk^0$.
For Calabi-Yau manifolds this reads as follows:

\begin{corollary}
--- Let $X$ be a Calabi-Yau manifold. The set $\kk^0$ of Ricci-flat
K\"ahler forms is linear if and only if for any Ricci-flat
K\"ahler form $\omega$ and any $\omega$-harmonic $(1,1)$-form
$\alpha$ the product $\alpha^2\omega^{N-2}$ is $\omega$-harmonic. In this
case, any class $\alpha$ in the positive cone $\kc_X\subset H^{1,1}(X,\IR)$ 
(cf.\ \cite{Huy}) is a K\"ahler class.\qed
\end{corollary}
\bigskip

{\bf Acknowledgement.} I wish to thank M.\ Lehn for a useful discussion
related to Sect.\ \ref{slgeneral}. The 
hospitality of the IHES, where the final version of this article
was prepared, is gratefully acknowledged.

\bigskip


\texttt{huybrech@mi.uni-koeln.de}
\end{document}